\documentclass[12pt, a4paper]{amsart}

\usepackage{fullpage}

\usepackage{color}
\definecolor{webcolor}{rgb}{0,0,1}
\definecolor{webbrown}{rgb}{.6,0,0}
\usepackage[
        colorlinks,
        linkcolor=webbrown, filecolor=webbrown,  citecolor=webbrown,
        pdfauthor={},
  pdftitle={},
]{hyperref}

\usepackage{type1cm}         
\usepackage{graphicx}        
\usepackage{multicol}        
\usepackage[bottom]{footmisc}

\usepackage{amsfonts, amsmath, amssymb}
\usepackage{xspace}          
\usepackage{ytableau}        
\usepackage{pgf}

\usepackage[draft=true]{minted} 
\setminted{breaklines=true}
\definecolor{lbcolor}{rgb}{0.9,0.9,0.9}
\setminted{bgcolor=lbcolor}
\setminted{fontsize=\footnotesize}

\usepackage{mathtools}
\usepackage{booktabs}        
\usepackage{tikz}            
\usepackage{csquotes}        

\usepackage[colorinlistoftodos, bordercolor=orange, backgroundcolor=orange!20, linecolor=orange, textsize=scriptsize]{todonotes} 

\usepackage[noend]{algorithmic}

\usepackage{calc,ifthen,alltt,tabularx,capt-of}
\usepackage[linesnumbered,longend,ruled,vlined]{algorithm2e}

\usepackage{longtable}
\usepackage{makecell}
\usepackage{wrapfig}

\usepackage[capitalise, noabbrev]{cleveref} 
\usepackage{enumitem} 

\newcommand\OSCAR{\texttt{OSCAR}\xspace}

\newcommand\FF{\mathbb{F}}
\newcommand\NN{\mathbb{N}} 
\newcommand\ZZ{\mathbb{Z}} 
\newcommand\QQ{\mathbb{Q}} 

\newcommand\RR{\mathbb{R}} 
\newcommand\CC{\mathbb{C}} 

\DeclareMathOperator{\Aut}{Aut} 
\DeclareMathOperator{\disc}{disc} 
\DeclareMathOperator{\Gal}{Gal} 
\DeclareMathOperator{\SL}{SL} 

\DeclareMathOperator{\Cl}{Cl}
\DeclareMathOperator{\TheProperSpec}{Spec}

\definecolor{promptColor}{rgb}{0.0,0.0,0.589}
\definecolor{brkpromptColor}{rgb}{0.589,0.0,0.0}
\definecolor{gapinputColor}{rgb}{0.589,0.0,0.0}
\definecolor{gapoutputColor}{rgb}{0.0,0.0,0.0}
\definecolor{darkgreen}{rgb}{0.05,0.6,0.1}
\definecolor{colrem}{rgb}{0,0.7,0}

\newcommand{\tmfloatcontents}{}
\newlength{\tmfloatwidth}
\newcommand{\tmfloat}[5]{
	\renewcommand{\tmfloatcontents}{#4}
	\setlength{\tmfloatwidth}{\widthof{\tmfloatcontents}+1in}
	\ifthenelse{\equal{#2}{small}}
	{\setlength{\tmfloatwidth}{0.45\linewidth}}
	{\setlength{\tmfloatwidth}{\linewidth}}
	\begin{minipage}[#1]{\tmfloatwidth}
		\begin{center}
			\tmfloatcontents
			\captionof{#3}{#5}
		\end{center}
\end{minipage}}

\newcommand{\fd}{.}

\title{Number Theory in Oscar}
\author{Claus Fieker}
\address{RPTU Kaiserslautern--Landau, Fachbereich Mathematik, Kaiserslautern, Germany.}
\email{claus.fieker@rptu.de}

\author{Tommy Hofmann}
\address{Universität Siegen, Naturwissenschaftlich-Technische Fakultät, Siegen, Germany.}
\email{tommy.hofmann@uni-siegen.de}


\usepackage{cite}
 
\begin{document}
\begin{abstract}
We give a brief introduction to computational algebraic number theory in
\OSCAR. Our main focus is on number fields, rings of integers and their invariants.
After recalling some classical results and their constructive counterparts,
we showcase the functionality in two examples related to the investigation of the Cohen--Lenstra heuristic for quadratic fields and
the Galois module structure of rings of integers.
\end{abstract}

\maketitle

\section{Introduction}
\label{ch:cs-number-theory}
The area of number theory is one of the oldest mathematical disciplines and
comprises different streams of research. Many developments have
been inspired by the study of Diophantine equations, which is the investigation of integer solutions
of polynomial equations with integer coefficients. One of the most famous examples is the
equation $x^n + y^n = z^n$ with $n \in \NN$, $n \geq 3$. The statement that this equation has no non-trivial integer
solution is commonly known as Fermat's last theorem, whose proof by Wiles~\cite{MR1333035} is
regarded as one of the pinnacles of the development of modern number theory in the 20th century.
Algebraic number theory, which is concerned with the investigation of number theoretic
questions using concepts from abstract algebra and counts Gauß~\cite{MR837656} among its
founding fathers, forms one of the pillars of modern number theory. The beauty
of this subject lies both in the depth of the
theory and the possibility to do concrete computations and
experiments.
In recent years, this fruitful interplay between theory and practical questions has been intensified further by the connection to the theory of secure communication, that is, cryptography.
In this chapter, we give a brief introduction to some of the
main characters from algebraic number theory and show how we can investigate them using the algorithmic tools provided by \OSCAR.

%

\label{sec:1}

\section{Number Fields and Rings of Integers}

One of the fundamental notions in algebraic number theory is that of a
\textit{number field}, which is a finite extension field of the rational numbers.
In this section we will explain the basics for working with these objects in \OSCAR (see also~\cite{MR3703682} on some background for the design choices of the number theory in \OSCAR).
We will freely use theoretical results from algebraic number theory, which can be found, for
example, in \cite{MR1282723, MR1697859, MR3822326}.
An algorithmic treatment of the subject can be found in \cite{MR1228206, MR1483321, MR1728313}.

\subsection{Construction of Number Fields and First Properties}

From theory we know that any number field $K$, that is, a finite extension of $\QQ$, can be generated
by a primitive element, say $K = \QQ(\alpha)$ for some $\alpha \in K$.
Since we are interested at this point only in the algebraic properties of $K$,
we can realize $K$ also as $K \cong \QQ[x]/(f)$, where $f \in \QQ[x]$ is the
minimal polynomial of $\alpha$ (which is necessarily irreducible).
Thus specifying a number field $K$ is as simple as specifying an irreducible rational polynomial, a \textit{defining polynomial} for $K$.
As an example, we want to construct the quadratic number field
$K = \QQ[x]/(x^2 - 235) \cong \QQ(\sqrt{5 \cdot 47})$.

\begin{minted}{jlcon}
julia> Qx, x = QQ["x"];

julia> K, a = number_field(x^2 - 235, "a")
(Number field of degree 2 over QQ, a)

julia> a^2 - 235 # confirm that a is a root of x^2 - 235
0
\end{minted}

Note that the construction of the number field returns as a second value a primitive element $\alpha$ satisfying $\alpha^2 - 235 = 0$. Since $\{1, \alpha\}$ is a $\QQ$-basis of $K$, we can now express arbitrary elements as linear combinations of $1$ and $\alpha$. Vice versa we can obtain the coordinates of any element.
Another way of constructing elements of $K$ is by mapping elements of $\QQ[x]$ under the canonical surjection $\QQ[x] \to \QQ[x]/(x^2 - 235) = K$.

\begin{minted}{jlcon}
julia> b = 2 + 1//3 * a
1//3*a + 2

julia> coordinates(b)
2-element Vector{QQFieldElem}:
 2
 1//3

julia> K(x^3 + x + 2) # map the polynomial x^3 + x + 2 to K
236*a + 2
\end{minted}

Of course, since a number field $K$ is foremost an extension of $\QQ$, we can ask for basic field theoretic properties of this extension as well as of the elements of $K$.

\begin{minted}{jlcon}
julia> degree(K)
2

julia> trace(a)
0

julia> norm(a)
-235
\end{minted}

\subsection{Three Theorems from Algebraic Number Theory}

For arithmetic purposes, one is often interested in the \textit{ring of integers} $\mathcal O_K$ of a number field $K$,
which is the subring comprised of all elements of $K$ with integral minimal polynomial.
Equivalently, $\mathcal O_K$ is the integral closure of $\ZZ$ in $K$. Because of this maximality property, it is also known as the \textit{maximal order} of $K$.
The structure of $\mathcal O_K$ is described by three classical theorems from algebraic number theory.
Before stating these theorems, we first have to introduce some notation.
The set $I_K$ of invertible fractional ideals of $\mathcal O_K$ forms an abelian group with respect to multiplication. It contains the subgroup $P_K$ of principal fractional ideals and the quotient $\Cl(K) = I_K/P_K$ is the \textit{class group} of $K$. The order $h_K = \lvert \Cl(K) \rvert$ is the \textit{class number} of $K$.
The \textit{signature} of $K$ is the tuple $(r, s)$ of non-negative integers, where $r$ denotes the number of real embeddings $K \to \RR$ and $2s$ the number of non-real embedding $K \to \CC$.
With this at hand, the ring, ideal and unit structure of $\mathcal O_K$ can be described as follows:

\begin{enumerate}
  \item\label{three:1}
    As an abelian group, $\mathcal O_K$ is a free $\ZZ$-module of rank $d = [ K : \QQ ]$. Thus there exists an \textit{integral basis} $(\omega_i)_{1 \leq i \leq d}$ satisfying $\mathcal O_K = \bigoplus_{i=1}^d \ZZ \omega_i$.
  \item\label{three:2}
    The ring $\mathcal O_K$ is a Dedekind domain with finite class group $\Cl(K)$.
  \item\label{three:3}
    The unit group $\mathcal O_K^\times$ is finitely generated of rank $r + s - 1$ with finite cyclic torsion subgroup, where $(r, s)$ is the signature of $K$. Thus we have $\mathcal O_K^\times \cong \ZZ^{r + s - 1} \times \ZZ/m\ZZ$, where $m \in \ZZ_{\geq 1}$ is the number of roots of unity of $K$.
\end{enumerate}

Since we approach this subject from a constructive perspective, these statements beg the questions:
Can we compute an integral basis, the class group or generators of the unit group?
We will now discuss how these questions can be answered with the help of \OSCAR.

\subsubsection*{The Ring of Integers}
\noindent
For statement~\ref{three:1}, note that already asking for the ring of integers of $K$ in \OSCAR will implicitly compute an integral basis, which can be queried for afterwards.

\begin{minted}{jlcon}
julia> OK = ring_of_integers(K);

julia> basis(OK)
2-element Vector{AbsSimpleNumFieldOrderElem}:
 1
 a

julia> discriminant(OK)
940
\end{minted}

Continuing with our example $K = \QQ[x]/(x^2 - 235)$, we get as a result $\{1, \alpha\}$, which means that
$\mathcal O_K = \ZZ \cdot 1 \oplus \ZZ\cdot \alpha$. From this we can also derive the discriminant $\operatorname{disc}(K)$ of $K$, which is equal to $\det((\operatorname{Tr}_{K/\QQ}(\omega_i \omega_j))_{1 \leq i, j \leq d}))$ and a measure for the ``size'' of the ring of integers.

We want to mention that this result has a simple geometric interpretation in case $K$ is an imaginary quadratic field.
In this case we can identify $K$ with a subset of $\CC$, such that the set $\mathcal{O}_K$ becomes a lattice (a free $\ZZ$-submodule of rank $2$). For example, if $L = \QQ(\alpha)$, $\alpha^2 + 3 = 0$, then $\mathcal{O}_L = \ZZ[(1 + \alpha)/2]$ has integral basis $(\omega_1, \omega_2) = (1, (1 + \alpha)/2)$.
We illustrate this using the Julia package \texttt{Plots.jl}

\begin{minted}{jlcon}
using Plots;
w1, w2 = Complex(1),(sqrt(Complex(-3)) + 1)/2
pts = filter(z -> max(abs(imag(z)), abs(real(z))) <= 3,
             [ a*w1 + b*w2 for a in -5:5 for b in -5:5]);
x = range(0, real(w1) + real(w2), 100);
y1 = (x -> min(sqrt(3)*x, imag(w2))).(x);
y2 = (x -> max(0, sqrt(3)*x - sqrt(3))).(x);
p = scatter(real.(pts), imag.(pts), framestyle=:origin,
            legend = false, mc = :red)
plot!(p, x, y1, fillrange = y2, fillalpha = 0.25, color = :blue)
\end{minted}

\begin{figure}[h]\centering
  \includegraphics[scale=0.45]{\fd/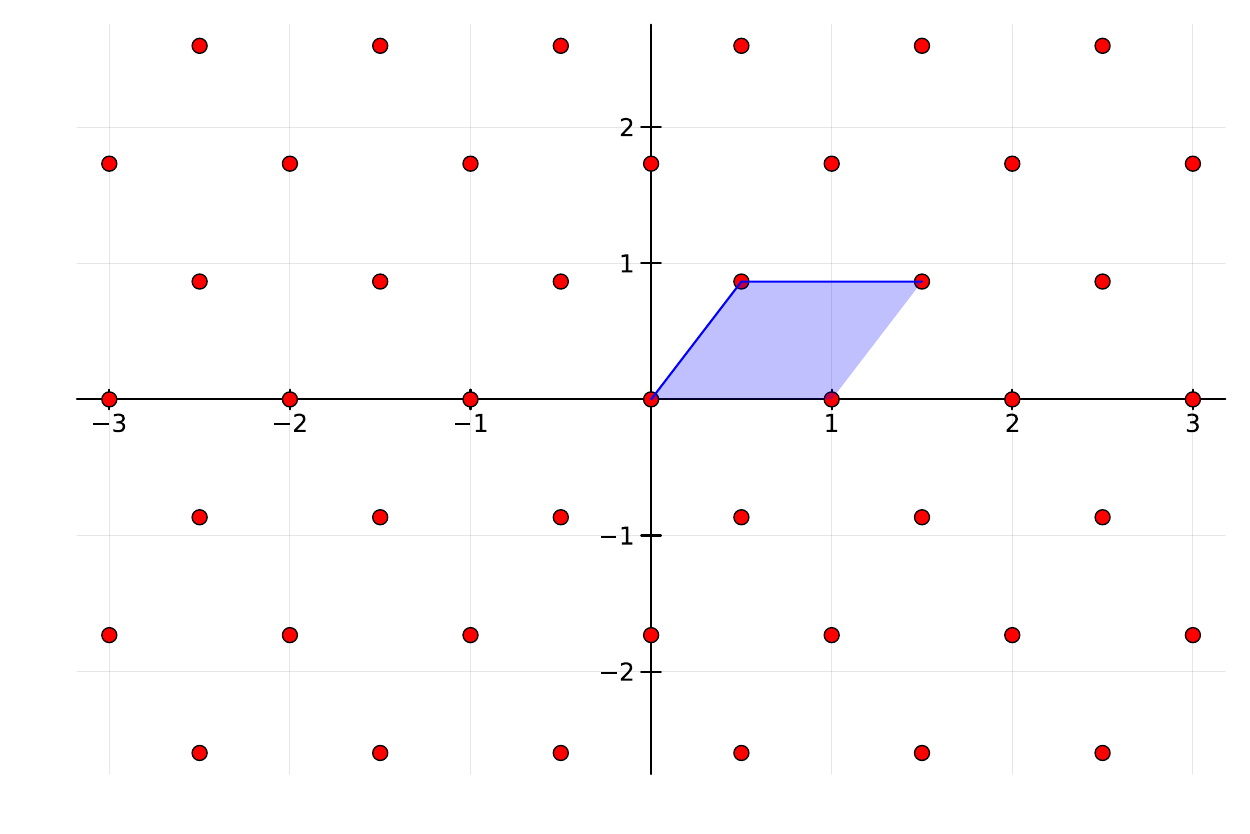}
  \caption{The ring $\ZZ[(1 + \sqrt{-3})/2]$ as a lattice of rank $2$ in $\CC$ and a fundamental domain.}
  \label{fig:oklattice}
\end{figure}

The result is displayed in Figure~\ref{fig:oklattice}. The red dots are the points of $\mathcal{O}_K$ with bounded real and imaginary part. The blue area is a fundamental domain for this lattice with area equal to $\sqrt{3}/2 = \sqrt{\lvert \operatorname{disc}(K) \rvert}/2$ (the factor of $2$ is introduced by an implicit choice of measure on $\CC$),
thus giving a precise meaning to the aforementioned claim on the relation between discriminant and ``size'' of the ring of integers.

\subsubsection*{Unique Factorization of Ideals and The Class Group}
\noindent
Let us now turn to statement~\ref{three:2}, which describes the ideal structure of the ring of integers.
That this ring is a Dedekind domain is equivalent to the unique factorization of non-trivial ideals
into prime ideals.
Before giving an example of this, we first consider the structure of prime ideals and highlight the celebrated theorem of Dedekind--Kummer, which, for almost all prime ideals, gives an explicit description in terms of properties of the defining polynomial of the number field.
Since $\mathcal O_K$ is an integral extension of $\ZZ$, the canonical map $\TheProperSpec(\mathcal O_K) \to \TheProperSpec(\ZZ)$ is surjective. Thus we obtain all prime ideals of $\mathcal O_K$ by determining the prime ideals lying over the rational primes in $\ZZ$.

\begin{minted}{jlcon}
julia> prime_ideals_over(OK, 7)
2-element Vector{AbsSimpleNumFieldOrderIdeal}:
 <7, a + 5>
 <7, a + 2>
\end{minted}

While this a priori only shows that $\mathfrak p_1 = (7, \alpha + 5)$ and $\mathfrak p_2 = (7, \alpha + 2)$ are the prime ideals lying above $7$, we have in fact the stronger statement that $7 \mathcal O_K = (7, \alpha + 4 ) \cdot (7, \alpha + 3)$. This is the factorization of $7 \mathcal O_K$ into prime ideals. We verify this by applying the theorem of Dedekind--Kummer,
which states that the prime ideal decomposition of $7 \mathcal O_K$ is in shape equal to the unique
factorization of $\overline{x^2 - 235} \in \FF_7[x]$ into irreducible polynomials, and that generators of the prime ideals are obtained by evaluating the lifted factors at $\alpha$.

\begin{minted}{jlcon}
julia> factor(change_coefficient_ring(GF(7), x^2 - 235))
1 * (x + 5) * (x + 2)
\end{minted}

The factorization into prime ideals is not limited to prime numbers. We can also factor arbitrary ideals, for example, $\alpha \mathcal O_K$.

\begin{minted}{jlcon}
julia> factor(a * OK)
Dict{AbsSimpleNumFieldOrderIdeal, Int64} with 2 entries:
  <47, a> => 1
  <5, a>  => 1
\end{minted}

We now turn our attention to the class group $\Cl(K)$. Let us first ask \OSCAR for this invariant and then see how this should be interpreted.

\begin{minted}{jlcon}
julia> A, m = class_group(OK);

julia> A
Z/6
\end{minted}

This tells us that the class group $\Cl(K)$ is, as an abelian group, isomorphic to the cyclic group $A \cong \ZZ/6\ZZ$.
If $A$ is not trivial, as in our case, this does not give a complete description of the group structure of $\Cl(K) = I_K/P_K$, since it does not tell us how the classes $[I], [J]$ of two (fractional) ideals $I, J$ of $\mathcal{O}_K$ behave under multiplication.
Here we encounter one of \OSCAR's design principles---the concept of \textit{interpretation maps}.
If the result of a computation carries the structure of, say, an abelian group, we will return an abstract abelian group together with an appropriate map (the \textit{interpretation map}), which is an isomorphism or can be used to reconstruct the isomorphism.
In our case, the returned map $m \colon A \to I_K$ is such that the composition
\[ A \xrightarrow[\phantom{xxx}]{m} I_K \longrightarrow I_K/P_K = \Cl(K) \]
with the natural projection is an isomorphism.
We can evaluate $m$ and also ask for preimages of arbitrary ideals, respectively, their classes.

\begin{minted}{jlcon}
julia> m(zero(A)) # apply m to the neutral element of A
<1, 1>
[...]

julia> P = prime_ideals_over(OK, 2)[1]
<2, a + 1>
[...]

julia> preimage(m, P)
Abelian group element [3]
\end{minted}

Thus we have shown that the ideal $\mathfrak p = (2, \alpha + 1)$ corresponds to the element $\overline 3 \in \ZZ/6\ZZ$ and thus has order $2$.
In particular this implies that $\mathfrak p^2 \in P_K$, that is, $\mathfrak p^2$ is principal.
The additional information about a possible generator of $\mathfrak p^2$ is not part of the map $m$, but
we can ask for it separately.

\begin{minted}{jlcon}
julia> is_principal_with_data(P^2)
(true, 2)

julia> P^2 == 2*OK
true
\end{minted}

\subsubsection*{The Unit Group}
\noindent
Finally we turn to statement~\ref{three:3}, which describes the structure of the unit group $\mathcal O_K^\times$. First note that if the signature of $K = \QQ[x]/(f)$ is $(r, s)$,
then the rational polynomial $f$ has $r$ real roots and $2s$ non-real complex roots. Hence the rank of $\mathcal O_K^\times$, which is equal to $r + s - 1$,
can be computed very easily. But as in the case of class groups, we also want to know an explicit isomorphism describing the structure of $\mathcal O_K^\times$ as an abelian group.

\begin{minted}{jlcon}
julia> U, mU = unit_group(OK);

julia> U
Z/2 x Z
\end{minted}

Thus $\mathcal O^\times_K \cong U = \ZZ \times \ZZ/2\ZZ$ as expected, since here $(r, s) = (2, 0)$. The interpretation map $U \to \mathcal O_K^\times$ is in this case an isomorphism. We can use it to determine explicit generators of the unit group, as well as, to map units to $U$, thus expressing elements in terms of the generators.

\begin{minted}{jlcon}
julia> mU(U[1]), mU(U[2])
(-1, 3*a + 46)

julia> preimage(mU, -214841715*a - 3293461126)
Abelian group element [1, 5]

julia> -214841715*a - 3293461126 == (-1)^1 * (3a + 46)^5
true
\end{minted}


Finally we want to mention how the unit group can be interpreted geometrically.
For the sake of simplicity we consider the number field $K = \QQ(\alpha)$ with $\alpha^2 - 3 = 0$, whose ring of integers is equal to $\ZZ[\alpha] = \{ a + b \alpha \mid a, b \in \ZZ\}$.
Since this a quadratic extension, we can embed $K$ into $\RR^2$ via the $\QQ$-linear map
\[ \varphi \colon K \to \RR^2, \ a + b\alpha \mapsto (a, b). \]
Now for any number field, the unit group is equal to $\mathcal{O}_K^\times = \{ x \in \mathcal{O}_K \mid \operatorname{N}_{K/\QQ}(x) = \pm 1 \}$.
For the field under consideration this implies
\[ \varphi(\ZZ[\alpha]^\times) = \{ (a, b) \in \ZZ^2 \mid a^2 - 3b^2 = \pm 1\}. \]
Thus $\ZZ[\alpha]^\times$ corresponds to the points of $H = \{ (a, b) \in \RR^2 \mid a = \pm \sqrt{3b^2 \pm 1}\}$ with integer
coordinates. We can plot this as follows:

\begin{minted}{jlcon}
pts = [ (a, b) for a in -8:8 for b in -5:5]
x = range(-7.5, 7.5, 500)
y1 = (x -> x^2 - 1 < 0 ? missing : sqrt((x^2 - 1)/3)).(x)
y2 = (x -> sqrt((x^2 + 1)/3)).(x)
scatter(pts, framestyle=:origin, leg = false, mc = :blue, alpha = 0.25)
plot!(x, [y1 -y1 y2 -y2], color = :black)
units = [(1, 0), (-1, 0), (2, 1), (-2, 1), (-2, -1),
         (2, -1), (7, 4), (7, -4), (-7, 4), (-7, -4)];
scatter!(units, color = :red)
\end{minted}

\begin{figure}[h]
  \includegraphics[scale=0.5]{\fd/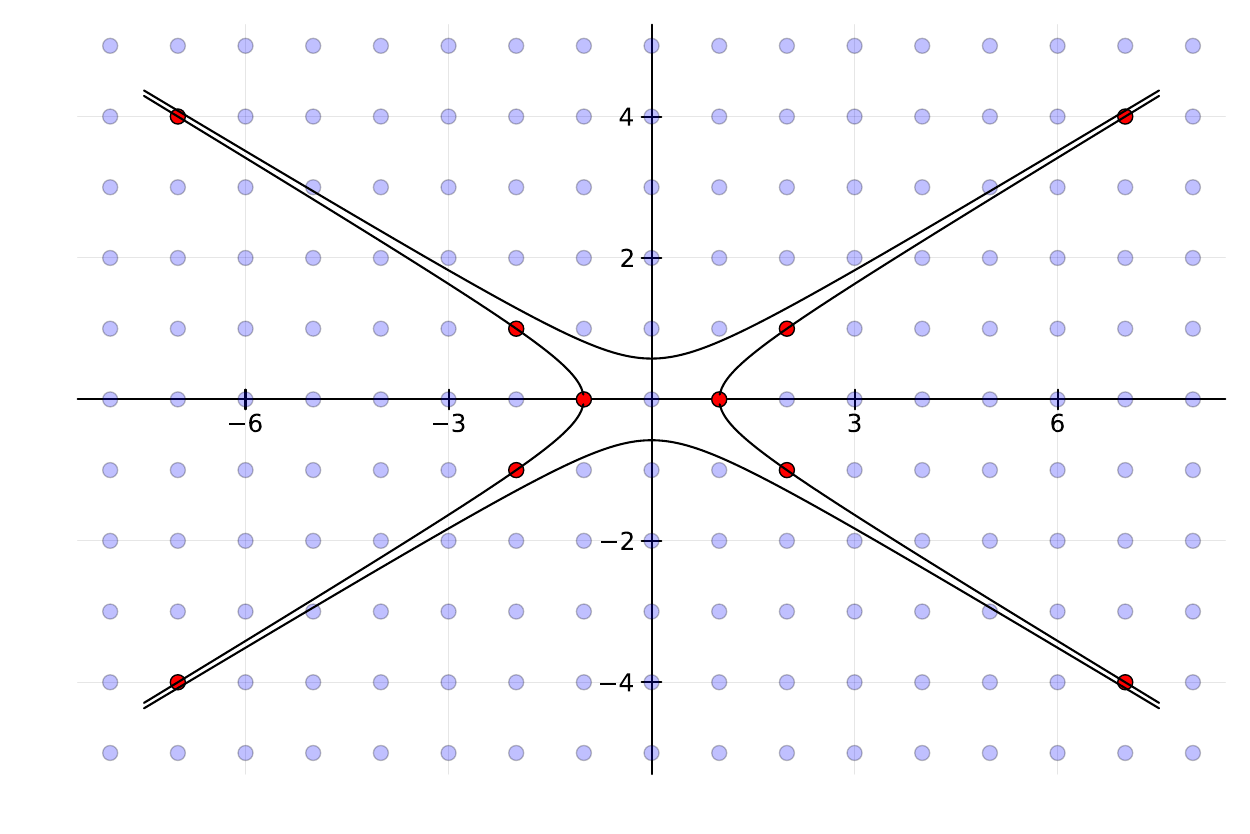}
  \caption{The unit group $\ZZ[\alpha]^\times$ interpreted as integer points lying on the union of two algebraic curves .}
  \label{fig:units}
\end{figure}

The result is depicted in Figure~\ref{fig:units}. The blue dots correspond to elements of $\ZZ[\alpha]$ with red dots corresponding to points lying in $H$, thus to units of $\ZZ[\alpha]$.
The points $(\pm 1, 0)$ correspond to the torsion units $\pm 1 \in \ZZ[\alpha]^\times$,
and $(\pm 2, \pm 1)$ correspond to the units $\pm 2 \pm \alpha$, which are in fact generators of the torsion-free part of $\ZZ[\alpha]^\times$.
The prominent symmetries of $\varphi(\ZZ[\alpha])$ (and of $H$) are explained by the existence of the two automorphisms $x \mapsto -x$ and $x \mapsto x^{-1}$ of $\ZZ[\alpha]^\times$.

Note that while Figure~\ref{fig:units} illustrates that units are quite likely to be ``sparse'' within the ring of integers, it does not depict the group structure of the unit group. To visualize the latter, we will make use
of the so-called logarithmic embedding (note that this historical name is slightly misleading, since this map is \textit{not} injective):
\[ \ell \colon K^\times \to \RR^2, \ a + b \alpha \mapsto (\log\lvert a + b \sqrt 3 \rvert, \log \lvert a - b \sqrt 3 \rvert). \]

The image $\ell(\ZZ[\alpha]^\times)$ is a free $\ZZ$-submodule of $\RR^2$ contained in the hyperplane $H = \{ (a, b) \in \RR^2 \mid a + b = 0 \}$, generated by $\ell(2 + \alpha) = (1.316..., -1.316...)$
We can plot this with the following commands:

\begin{minted}{jlcon}
l(a, b) = (log(abs(a + sqrt(3)*b)), log(abs(a - sqrt(3)*b)))
pts = filter!(z -> maximum(abs.(z)) < 3,
              [l(a, b) for a in -50:50 for b in -50:50])
scatter(pts, framestyle=:origin, leg = false, mc = :blue, alpha = 0.25)
x = range(-3, 3, 100); y = -x; plot!(x, y, color = :black)
units = [ i .* l(2, 1) for i in -2:2]
scatter!(units, color = :red)
\end{minted}

\begin{figure}[h]
  \includegraphics[scale=0.5]{\fd/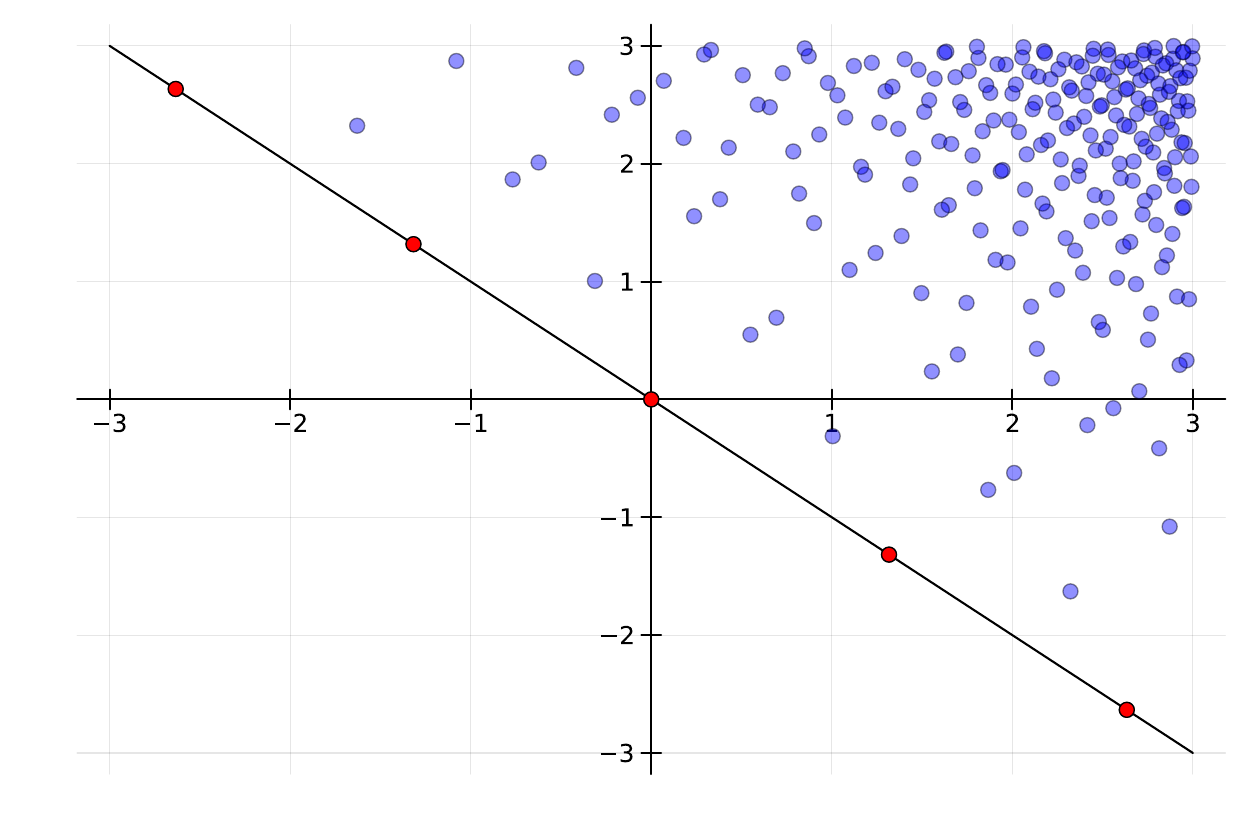}
\caption{The unit group $\ZZ[\alpha]^\times$ under the logarithmic embedding $\ell$.}
  \label{fig:logunits}
\end{figure}

The result is depicted in Figure~\ref{fig:logunits}. As before, the blue dots correspond to elements of $\ZZ[\alpha]$ with red dots corresponding to points lying on $H$, that is, corresponding to units of $\ZZ[\alpha]$.
The abelian structure of $\ZZ[\alpha]$ is now clearly visible. Note that, in comparison to Figure~\ref{fig:units}, we have ``lost'' one symmetry. This is explained by the fact that $-1 \in \ker(\ell)$, and hence the automorphism $x \mapsto -x$ of $\ZZ[\alpha]^\times$ becomes the identity.

\section{General Number Fields and Morphisms}\label{sec:genhom}

So far we have constructed number fields only as simple extensions of $\QQ$, where we have provided the minimal polynomial of a generator.
Although due to the primitive element theorem, any number field is of this kind,
in practice, number fields are not always given in this form. Often, in the context of the specific problem under consideration, there is a description as a non-simple extension of $\QQ$ or of other number fields.
In this regard, \OSCAR is quite flexible in that a number field is allowed to be either
\begin{itemize}
  \item
    a simple extension $\QQ(\alpha)$ or a non-simple extension $\QQ(\alpha_1,\dotsc,\alpha_n)$ of $\QQ$ (referred to as an \textit{absolute field}); or
  \item
    a simple extension $K(\alpha)$ or a non-simple extension $K(\alpha_1,\dotsc,\alpha_n)$, where $K$ is a previously defined number field (referred to as \textit{a relative field}).
\end{itemize}

Let us illustrate this by constructing the biquadratic non-simple extension
$K = \QQ(\sqrt 2, \sqrt 3)$. This is as simple as providing the defining
polynomials of the generators, which in our case are $x^2 - 2$ and $x^2 - 3$.

\begin{minted}{jlcon}
julia> K, a = number_field([x^2 - 2, x^2 - 3]);

julia> a[1]^2 == 2 && a[2]^2 == 3
true
\end{minted}

To illustrate the relative fields, we now construct $\QQ(\sqrt 2, \sqrt 3)$ as the quadratic extension $L = k(\sqrt 3)$, where $k = \QQ(\sqrt 2)$.
Note that although $t^2 - 2 \in k[t]$ has rational coefficients, it is constructed as a polynomial over $k$.
As a consequence, when we use it to define a number field, we obtain an extension of $k$.

\begin{minted}{jlcon}
julia> k, b = quadratic_field(3);

julia> kt, t = k["t"];

julia> L, c = number_field(t^2 - 2);
\end{minted}

Note that although $L$ and $K$ are isomorphic as extensions of $\QQ$, they behave very differently.
Since $L$ was created as an extension of $k$, all the field theoretic
properties of $L$ and its elements are relative to the base field $k$.

\begin{minted}{jlcon}
julia> base_field(L) == k
true

julia> base_field(K) == QQ
true

julia> degree(K)
4

julia> degree(L)
2
\end{minted}

We can still ask the system for information of the number field $L$, regarded as an absolute field, that is, as a field extension of $\QQ$.

\begin{minted}{jlcon}
julia> absolute_degree(K)
4

julia> absolute_norm(c)
4
\end{minted}

The flexible framework to construct number fields as arbitrary extensions of number fields is accompanied by
a framework for working with maps between number fields of similar flexibility.
Specifying a map $L \to M$ from a number field $L$ with base field $k$ is done by providing the image(s) of the generator(s)
and an optional map $k \to M$.
We illustrate this by constructing the map defined by $\sqrt 2 \mapsto - \sqrt 2$ and $\sqrt 3 \mapsto -\sqrt 3$ for both fields $K$ and $L$ encountered before.

\begin{minted}{jlcon}
julia> f = hom(K, K, [-a[1], -a[2]]);

julia> f(a[1]) == -a[1] && f(a[2]) == -a[2]
true

julia> g = hom(L, L, hom(k, k, -b), -c);

julia> g(c) == -c && g(L(b)) == L(-b)
true
\end{minted}


As we have seen, the construction of number fields is quite flexible and allows for arbitrary towers of number fields.
On the other hand, it is often necessary to collapse a given number field $L$ into an absolute simple extension, that is, into the form $\QQ(\alpha)$ for some $\alpha$.
This need arises, for example, when a certain functionality is provided only for absolute simple number fields.
We illustrate this by computing the class group of $L = \QQ(\sqrt 3)(\sqrt 2)$ and proving that $h_L = 1$.

\begin{minted}{jlcon}
julia> Lprime, LprimetoL = absolute_simple_field(L);

julia> defining_polynomial(Lprime)
x^4 - 10*x^2 + 1
\end{minted}

The first return value  $L'$ is a simple extension of $\QQ$ and the second return value is a $\QQ$-linear isomorphism $L' \to L$. We now determine the class number of $L$ via $L'$:

\begin{minted}{jlcon}
julia> C, _ = class_group(ring_of_integers(Lprime)); order(C)
1
\end{minted}

Note that similar functionality exists to collapse a non-simple extension into a simple extension (\texttt{simple\_extension}) and to
collapse a relative extension $L/K/k$ into an extension $L'/k$ (\texttt{collapse\_top\_layer}).

\section{Embeddings of Number Fields}

We now want to stress an important point about computational aspects of number fields, which requires some adjustment
to the conventional way of working with number fields. If we define the
number field $K$ as a field extension $\QQ[x]/(f)$ of $\QQ$ for some irreducible polynomial $f \in \QQ[x]$,
then we do not consider $K$ to be embedded into $\CC$. 
For example, if $K = \QQ[x]/(x^2 - 2)$, then the primitive element $\alpha = \overline x
\in K$ satisfies $\alpha^2 - 2 = 0$, but we cannot tell whether $\alpha$ is equal to $\sqrt
2$ or $-\sqrt 2$.
On the theoretical side, this ambiguity comes from the fact that for a number field $K$, there is in general more than one morphism $K \to \CC$.

\subsection{Embeddings}

From the standpoint of embeddings and places, the presence of more than one
morphism $K \to \RR$ is seen as additional arithmetic information of the number
field.
A \textit{complex embedding} of a number field $K$ is a morphism $\iota \colon K \to \CC$.
If $\iota(K) \subseteq \RR$, the embedding is called \textit{real}, and \textit{imaginary} otherwise.
If $K = \QQ(\alpha) = \QQ[x]/(f)$ is a simple absolute extension, imaginary embeddings correspond to the roots of $f$ in $\CC \setminus \RR$, with the real embeddings corresponding to real roots.

\begin{minted}{jlcon}
julia> K, a = number_field(x^3 - 2, "a");

julia> signature(K)
(1, 1)

julia> real_embeddings(K)
1-element Vector{AbsSimpleNumFieldEmbedding}:
 Complex embedding corresponding to 1.26 of K

julia> embs = complex_embeddings(K) # printed with limited precision
3-element Vector{AbsSimpleNumFieldEmbedding}:
 Complex embedding corresponding to 1.26 of K
 Complex embedding corresponding to -0.63 + 1.09 * i of K
 Complex embedding corresponding to -0.63 - 1.09 * i of K

julia> embs[1](a - 1)
[0.2599210499 +/- 8.44e-11]
\end{minted}

\subsection{Embedded Number Fields}

We continue with the number field $K = \QQ[x]/(x^3 - 2)$ and primitive element $\alpha$ satisfying $\alpha^3 - 2 = 0$.
To obtain an ordered field, where $\alpha$ is in fact equal to $\sqrt[3]{2}$, we can construct an embedded number field, where we additionally supply an approximation of the root during the construction.

\begin{minted}{jlcon}
julia> K, a = embedded_number_field(x^3 - 2, 1.26);

julia> 0 < a < 2
true
\end{minted}

Note that such an embedded number field does not support number theoretic operations, but is useful when there is a need for an ordered field, for example, in polyhedral geometry.

\section{Symmetries}

The symmetries of a number field arise as a purely field-theoretic concept and play an fundamental role in
algebraic number theory.
On the one hand, these symmetries are an important tool when studying arithmetic objects connected to or defined over number fields.
On the other hand, various deep open problems in algebraic number theory are related to the study of the symmetries of number fields themselves, one of the most prominent being the Inverse Galois Problem.
In this section we explain how to compute and work with symmetries of number fields, which come in two different flavors: automorphism groups and Galois groups.

\subsection{Automorphisms}

Let $K/\QQ$ be a number field. The automorphism group $\Aut_\QQ(K)$ is by definition equal to the set of all ring isomorphisms $K \to K$.
Let us illustrate this using the number field $K/\QQ$ with defining polynomial
$x^4 - 13x^2 + 16$.

We can easily compute this with \OSCAR as follows.
After defining the field, we tell the system to compute the automorphism group as a permutation group.

\begin{minted}{jlcon}
julia> K, a = number_field(x^4 - 13*x^2 + 16, "a");

julia> G, m = automorphism_group(PermGroup, K); G # we only print G
Permutation group of degree 4
\end{minted}

(Note that the abstract permutation group $G$ representing the automorphism group is not unique. Due to the randomized nature of the algorithm, the permutations generating the group can differ.)
Since the two permutations commute and are of order $2$, we see that $G$ is isomorphic to $\mathrm{C}_2 \times \mathrm{C}_2$, which we can also confirm as follows:

\begin{minted}{jlcon}
julia> describe(G)
"C2 x C2"
\end{minted}

So far we have not seen any automorphism of $K$, but only some abstract description of the group of all automorphisms.
This additional information is provided by the map $m$, which is an isomorphism $G \to \Aut_\QQ(K)$.

\begin{minted}{jlcon}
julia> m(G[1])
Map
  from number field of degree 4 over QQ
  to number field of degree 4 over QQ

julia> m(G[1])(a)
-1//4*a^3 + 13//4*a
\end{minted}

These automorphisms are just special cases of morphisms that we already encountered in Section~\ref{sec:genhom}. For example, we can construct the automorphism of $K$ defined by $a \mapsto -a$ and ask for its preimage in $G$:

\begin{minted}{jlcon}
julia> preimage(m, hom(K, K, -a))
(1,2)(3,4)
\end{minted}

In agreement with the overall design principles for number fields, the automorphism group respects also the (possible) relative structure of a number field. For example, consider $k = \QQ(\alpha)$ with defining polynomial $x^2 - 18$ and $K = k(\beta)$ with defining polynomial $x^4 + (\alpha + 6)x^2 + 2\alpha + 9$.


\begin{minted}{jlcon}
julia> k, a = number_field(x^2 - 18, "a"); kt, t = k["t"];

julia> K, b = number_field(t^4 + (a + 6)*t^2 + 2a + 9, "b");

julia> G, m = automorphism_group(PermGroup, K); describe(G)
"C4"
\end{minted}

Thus $\Aut_k(K) \cong G \cong \mathrm{C}_4$ is cyclic of order $4$.
While any $k$-linear automorphism is necessarily $\QQ$-linear, not every $\QQ$-linear automorphism must be $k$-linear, and in fact, it might not even stabilize $k$. Thus, in general we have $\Aut_k(K) \subsetneq \Aut_{\QQ}(K)$.
As we will see now, the field $K$ has this property.
The group $\Aut_{\QQ}(K)$, which we refer to as the \textit{absolute automorphism group} can be computed as follows.

\begin{minted}{jlcon}
julia> G, m = absolute_automorphism_group(PermGroup, K); describe(G)
"Q8"
\end{minted}


\subsection{Galois Groups}

The second kind of symmetry we will discuss are Galois groups of number fields.
We will illustrate this using the number field $K = \QQ(\alpha)$ with defining polynomial $x^4 - 2$.

\begin{minted}{jlcon}
julia> K, a = number_field(x^4 - 2);

julia> G, mA = automorphism_group(PermGroup, K);

julia> describe(G)
"C2"
\end{minted}

Given that this field is not normal, the automorphism group of $K$ is not ``large enough''.
Here we only have one non-trivial automorphism, sending $\alpha$ to $-\alpha$.
One way to overcome this apparent lack of automorphisms is to pass to the Galois group, which is by
definition the automorphism group of a splitting field of a defining polynomial $f$ of $K$.
In \OSCAR we can directly ask for the Galois group.

\begin{minted}{jlcon}
julia> G, C = galois_group(K); describe(G)
"D8"
\end{minted}

We can verify that this is (isomorphic to) the automorphism group of a splitting field of $x^4 - 2$.

\begin{minted}{jlcon}
julia> L = splitting_field(x^4 - 2);

julia> GL, = automorphism_group(PermGroup, L);

julia> fl, = is_isomorphic(G, GL); fl
true
\end{minted}

The automorphism group of the splitting field of $f$ acts naturally on the $d = \deg(f)$ roots of $f$ and is therefore
naturally isomorphic to a subgroup of the symmetric group $\mathfrak{S}_d$ on $d$ letters.
For this reason the permutation group $G$ computed is not any permutation group isomorphic to the Galois group,
but a group of explicit permutations of the roots of the defining polynomial in a suitable field.
The information about the roots can be accessed via the second return value of the \texttt{galois\_group} command.

\begin{minted}{jlcon}
julia> rts = roots(C, 2)
4-element Vector{QadicFieldElem}:
 (8*11^0 + O(11^2))*a + 8*11^0 + 9*11^1 + O(11^2)
 (3*11^0 + 10*11^1 + O(11^2))*a + 7*11^0 + 4*11^1 + O(11^2)
 (3*11^0 + 10*11^1 + O(11^2))*a + 3*11^0 + 11^1 + O(11^2)
 (8*11^0 + O(11^2))*a + 4*11^0 + 6*11^1 + O(11^2)

julia> parent(rts[1])
Unramified extension of 11-adic numbers of degree 2
\end{minted}

We can read off that the roots of $x^4 - 2$ are approximated in a local field, more precisely, the unramified quadratic extension of the field of $11$-adic numbers $\QQ_{11}$.
The approximated roots can be computed with arbitrary precision---for the sake of simplicity we only asked for a precision of $2$.
An interesting side-effect of Stauduhar's algorithm~\cite{MR327712} (see also~\cite{MR4490928} for details on the implementation in \OSCAR), which is the algorithm underlying the Galois group computation, is that it allows us to determine subfields of the splitting field exactly, still using only approximated roots.
As an illustration we consider the field $K$ with defining polynomial $x^9 - 3x^8 + x^6 + 15x^5 - 13x^4 - 3x^3 + 4x - 1$. The Galois group of this field is the group $G = \operatorname{ASL}_2(\FF_3)$, the affine special linear group on $\FF_3^2$ of order $216$. This group has a subgroup $H$ isomorphic to the Heisenberg group $\mathrm{He}_3$ of order $27$ with $G/\operatorname{Core}(H) \cong \SL_2(\FF_3)$ of order $24$.
Thus the fixed field $k = K^H$ is non-normal of degree $8$ with Galois group $\SL_2(\FF_3)$.
We can compute a number field (and therefore a defining polynomial) with this Galois group as follows:

\begin{minted}{jlcon}
K, a = number_field(x^9 - 3*x^8 + x^6 + 15*x^5 - 13*x^4 -
                    3*x^3 + 4*x - 1, "a")
G, C = galois_group(K)
subsG = subgroups(G)
H = first(H for H in subsG if order(H) == 27)
k, = simplify(fixed_field(C, H))
\end{minted}

\begin{minted}{jlcon}
julia> defining_polynomial(k)
x^8 + x^7 - 7*x^6 - 23*x^5 + 22*x^4 + 80*x^3 + 99*x^2 + 27*x + 8

julia> describe(galois_group(k)[1])
"SL(2,3)"
\end{minted}

Thus we have found the field $k$. Note that the same computation using automorphism groups is
prohibitively expensive, since it would require the explicit construction of a normal closure of $K$, which has degree $216$.

We end this section by highlighting a few of the algorithmic differences between Galois group and automorphism group computations of a normal number field $K$ with defining polynomial $f \in \QQ[x]$.
Although in this situation
the automorphism group is isomorphic
to the Galois group, the Galois group is still given as a permutation group acting on the (approximated) roots in a local splitting field, while the automorphism
group will implicitly use the roots of $f$ in $K$. As a consequence, the actual size of the Galois group is of much less importance for Galois group computations, while being critical for automorphism group computations.
On the other hand, automorphisms can be applied directly to elements and ideals, while elements of
the Galois group actually do not act on the field.

In comparison to most older implementations, the Galois group computation in
\OSCAR does not use any tables of precomputed data, hence, subject to
runtime and memory constrains, can be applied to fields and polynomials
of arbitrary degree, as well as to reducible polynomials.



\section{Examples}

We end this chapter with some \OSCAR applications to a few well-known problems
from algebraic number theory.

\subsection{Class Group Heuristics}

One of the striking applications of computational tools for algebraic number theory
has been experiments with arithmetic properties in families of number fields or related objects.
The tabulation of inequivalent binary quadratic forms in certain discriminant ranges carried out by Gauß in article 303 of his ``Disquisitiones'' in 1801, which led him to formulate various conjectures about class numbers, is one of the earliest works in this direction.

Here we consider the distribution of class groups in families of quadratic fields,
which goes back to seminal work of Cohen--Lenstra~\cite{MR756082}.
More precisely,
for $X \in \RR$, $X > 0$, we will consider the family
\[ \mathcal{F}(X) := \{K/\QQ \text{ real quadratic}, \ \lvert \disc(K) \rvert \leq X\} \]
of real quadratic fields with bounded discriminant inside a fixed algebraic closure of $\QQ$.
If $p$ is an odd prime and $A$ a finite $p$-group, then the Cohen--Lenstra heuristic predicts that
\[ \lim_{X \to \infty} \frac{\lvert \{K \in \mathcal{F}(X) \mid \operatorname{Cl}(K)[p^\infty] \cong A \}\rvert}{\lvert \mathcal{F}(X) \rvert} = \frac{w_p}{\lvert A \rvert \cdot \lvert \Aut(A)\rvert}, \]
where $w_p = \prod_{k \geq 2} (1 - p^{-k})$ and $\operatorname{Cl}(K)[p^\infty]$ is the Sylow $p$-subgroup of the class group.
One can interpret the left hand side as the probability that the Sylow $p$-subgroup of the class group of a random real quadratic field is isomorphic to $A$.
Our aim is to numerically investigate this conjecture for $p = 5$ and $X = 10^6$.
(Some of the computations described below require a few minutes of runtime. At the expense of working with a smaller data set, one can choose $X = 10^{5}$ to obtain weaker results more quickly.)
We begin by computing the set of quadratic fields $\mathcal{F}(X)$ and their class groups.

\begin{minted}{jlcon}
julia> fields = [d for d in 1:10^6 if
                 d != 1 && is_fundamental_discriminant(d)];

julia> length(fields)
303957

julia> clgrps = [class_group(quadratic_field(d, cached = false)[1])[1]
                 for d in fields];
\end{minted}

To obtain the statistics, we will make use of the Julia package \texttt{Tally.jl}:

\begin{minted}{jlcon}
julia> using Tally;

julia> tally(clgrps,
             by = x -> sylow_subgroup(x, 5)[1],
             equivalence = (x, y) -> is_isomorphic(x, y)[1])
Tally with 303957 items in 4 groups:
[Z/1]     | 291400 | 95.869%
[Z/5]     |  12324 |  4.055%
[Z/25]    |    230 |  0.076%
[(Z/5)^2] |      3 |  0.001%
\end{minted}

We now determine the proportions predicted by the conjecture. We first construct the abelian groups (by specifying elementary divisors) and then determine the values.

\begin{minted}{jlcon}
julia> ab_grps = abelian_group.([[1], [5], [25], [5, 5]]);

julia> w = prod(1 - 1/5.0^i for i in 2:1000)
0.9504159948390403

julia> [ 100 * w/(1.0 * (order(A) * 
         order(automorphism_group(A)))) for A in ab_grps]
4-element Vector{BigFloat}:
 95.041[...]
  4.752[...]
  0.190[...]
  0.007[...]
\end{minted}

Thus we see that using the first 300\,000 fields, the numbers
match approximately.
We now want to have a closer look at how the discrepancy between the predicted and actual value evolves with increasing $X$.
For the sake of simplicity, we consider the probability that the class number of a real quadratic field is divisible by $5$. The Cohen--Lenstra heuristic predicts this value to be
\[ \lim_{X \to \infty} \frac{\lvert \{K \in \mathcal{F}(X) \mid 5 \text{ divides } h_K \}\rvert}{\lvert \mathcal{F}(X) \rvert} = 1 - w_5 = 1 - \prod_{k=2}^\infty \left(1 - \frac{1}{5^k}\right). \]

We now collect the predicted and actual proportions for various values of $X$.

\begin{minted}{jlcon}
julia> clnumbs = order.(clgrps);

julia> x = range(1, length(clnumbs), 1000);

julia> values = [count(c -> is_divisible_by(c, 5),
                       clnumbs[1:Int(ceil(step))])/Int(ceil(step)) for
                 step in x];
\end{minted}

\begin{figure}[h]\centering
  \includegraphics[scale=0.5]{\fd/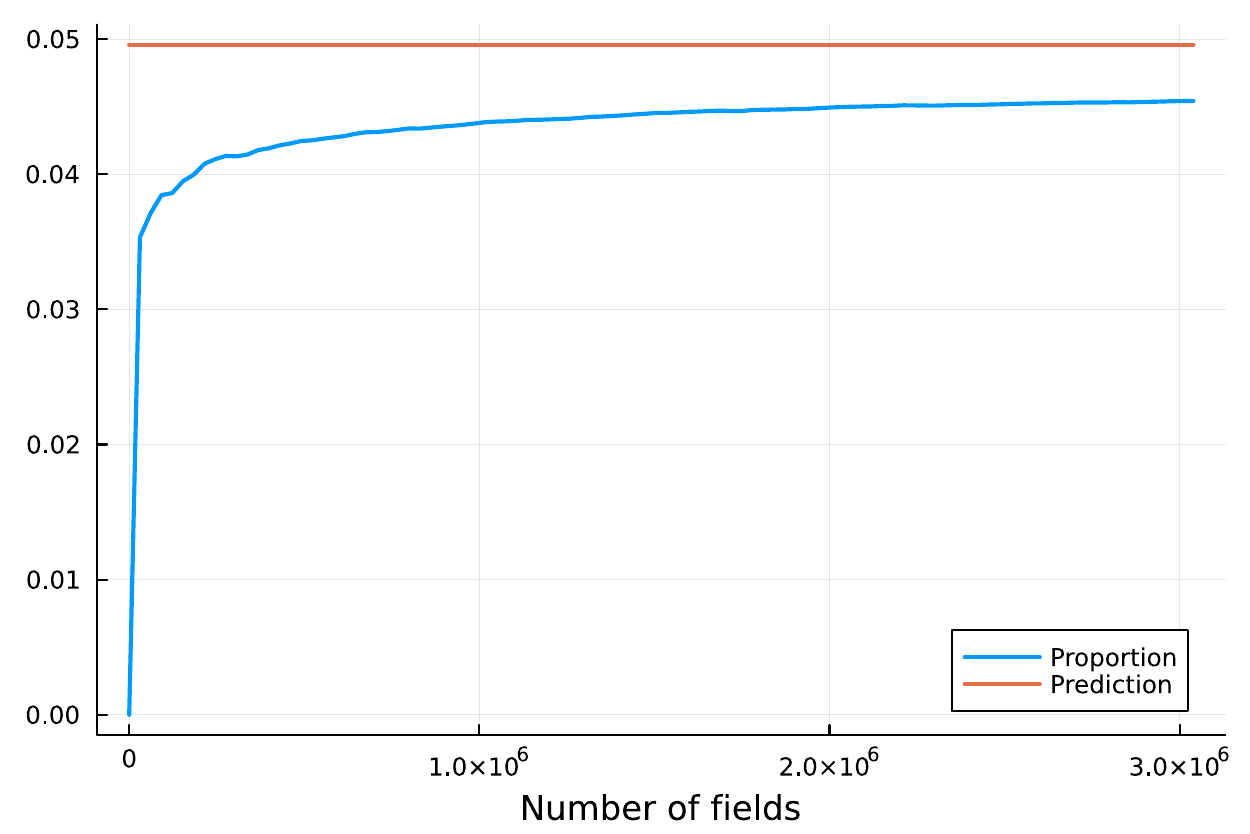}
  \caption{Predicted and actual proportion of real quadratic fields with class number divisible by $5$ with discriminant bounded by $10^7$.}
  \label{fig:cl}
\end{figure}

The plot of this result is presented in Figure~\ref{fig:cl} for $X \leq 10^7$ (for this bound, the computation of the class groups takes a day).
Note that we plot the proportion against $\lvert \mathcal F(X) \rvert$ instead of $X$.

\begin{minted}{jlcon}
julia> pr = 1 - prod(1 - 1/5.0^i for i in 2:1000)
0.04958400516095973

julia> using Plots;
       plot(x, [values fill(pr, length(x))], linewidth = 2,
            xlabel = "Number of fields",
            label = ["Proportion" "Prediction"])
\end{minted}

Finally we note that the Cohen--Lenstra heuristic has been extended
to arbitrary number fields and relative extensions by Cohen--Martinet~\cite{MR866103} with modifications due to Malle~\cite{MR2441080, MR2778658}; see also~\cite{Sawin2023} for some recent progress.

%

\subsection{Galois Module Structure}

For a normal extension $K/\QQ$, the Galois group $G = \Gal(K/\QQ)$ (which we interpret simply as the automorphism group $\Aut_\QQ(K)$)
acts naturally on many arithmetic objects derived from or defined over $K$.
The theory of Galois modules is an area of number theory concerned with the study of such actions.

\subsubsection*{The Normal Basis Theorem}
\noindent
The simplest example of a Galois module is the additive group of the field $K$ itself, on which $G$ acts naturally as automorphism group.
In this case the question about the Galois module structure is answered
by the celebrated normal basis theorem, which asserts that there exists an element $\beta \in K$
such that $(\sigma(\beta))_{\sigma \in G}$ is a basis of $K$ when considered as a $\QQ$-vector space. We refer to such an element $\beta$ as a
\textit{normal basis generator}.
There are different ways to find such an element for a given number field. We present a simple randomized approach and illustrate
it with $K = \QQ(\alpha)$, where $\alpha^4 + 4\alpha^2 + 2 = 0$.
Note that $K/\QQ$ is a Galois extension.
Since there are abundantly many normal basis generators, the idea is to pick a random element $\beta \in K$ and verify the condition of the definition:

\begin{minted}{jlcon}
julia> K, a = number_field(x^4 + 4*x^2 + 2, "a");

julia> b = rand(K, -10:10) # random element with coefficients
                           # between -10 and 10
5*a^3 - 7*a^2 + 8*a - 8

julia> A, mA = automorphism_group(K);
\end{minted}

Now construct the list $(\sigma(\beta))_{\sigma \in G}$:

\begin{minted}{jlcon}
julia> list = [mA(s)(b) for s in A]
4-element Vector{AbsSimpleNumFieldElem}:
 5*a^3 - 7*a^2 + 8*a - 8
 -7*a^3 + 7*a^2 - 26*a + 20
 -5*a^3 - 7*a^2 - 8*a - 8
 7*a^3 + 7*a^2 + 26*a + 20
\end{minted}

To verify that these elements are $\QQ$-linearly independent, we consider the coordinates with respect to the canonical basis:

\begin{minted}{jlcon}
julia> X = matrix(QQ, coordinates.(list))
[-8     8   -7    5]
[20   -26    7   -7]
[-8    -8   -7   -5]
[20    26    7    7]
\end{minted}

Thus, we have $X \cdot (1, \alpha, \alpha^2, \alpha^3)^t = (\sigma_1(\beta), \sigma_2(\beta), \sigma_3(\beta), \sigma_4(\beta))^t$ and it remains to check whether $X$ is invertible.

\begin{minted}{jlcon}
julia> det(X) != 0
true
\end{minted}

This proves that $\beta =  -2\alpha^3 + 10\alpha^2 - 2\alpha$ is indeed a normal basis generator of $K$. Note that not every primitive element is a normal basis generator. For this given quartic field $K$, the element $\alpha$ is, for example, not a normal basis generator:

\begin{minted}{jlcon}
julia> det(matrix(QQ, [coordinates(mA(s)(a)) for s in A]))
0
\end{minted}

Finally we mention that the normal basis theorem has an important module theoretic interpretation. It asserts that the number field $K$ is a free $\QQ[G]$-module of rank $1$ and any element $\beta \in K$ with $K = \QQ[G] \cdot \beta$ is a normal basis generator.

\subsubsection*{Integral Normal Bases}
\noindent
Using the same setup as before, we now consider the integral structure of a normal number field $K/\QQ$.
Since integrality is preserved under automorphisms, the Galois group $G = \Gal(K/\QQ)$ acts on $\mathcal{O}_K$, turning the ring of integers $\mathcal{O}_K$ into a $\ZZ[G]$-module.
In view of the normal basis theorem it is only natural to ask whether $\mathcal{O}_K$ is a free $\ZZ[G]$-module (necessarily of rank $1$).
An affirmative answer is equivalent to the existence of a \textit{normal integral basis generator}, that is, an element $\beta \in \mathcal{O}_K$ such that $(\sigma(\beta))_{\sigma \in G}$ is a $\ZZ$-basis of $\mathcal{O}_K$.
While it is in general a hard problem to decide freeness of $\ZZ[G]$-modules, \OSCAR provides some functionality to work with such objects.
We want to investigate whether the quartic number field $K = \QQ(\alpha)$ defined before has a normal integral basis generator.

To this end, we first construct a $\QQ[G]$-linear isomorphism $f\colon K \to V$, where $V$ is the regular $\QQ[G]$-module (the existence of which is guaranteed by the normal basis theorem)

\begin{minted}{jlcon}
julia> V, f = galois_module(K);
\end{minted}

Next, we translate the ring of integers to a $\ZZ[G]$-submodule $M$ of $V$ and ask whether it is free or not as a $\ZZ[G]$-module.

\begin{minted}{jlcon}
julia> OK = ring_of_integers(K);

julia> M = f(OK);

julia> is_free(M)
false
\end{minted}

Note that we could have already seen this by investigating whether $M$ is locally free at $2$, that is, whether $M \otimes \ZZ_2$ is a free $\ZZ_2[G]$-module.

\begin{minted}{jlcon}
julia> fl = is_locally_free(M, 2)
false
\end{minted}

In fact, the local non-freeness is equivalent to an arithmetic property of $\mathcal{O}_K$. Noether's criterion asserts that $\mathcal{O}_K$ is locally free at all primes $p$ if and only if $K/\QQ$ is tamely ramified.
The latter condition means that for all ramified primes $p$ and prime ideals $\mathfrak{p}$ lying above $p$, the ramification index $e(\mathfrak p/p)$ is not divisible by $p$. This is a criterion that can easily be checked:

\begin{minted}{jlcon}
julia> prime_decomposition_type(OK, 2)
1-element Vector{Tuple{Int64, Int64}}:
 (1, 4)
\end{minted}

Thus $2 \cdot \mathcal{O}_K = \mathfrak p_2^4$ for some prime ideal $\mathfrak p_2$, $2 \mid 4 = e(\mathfrak p_2/2)$ and the extension is therefore not tamely ramified.
If we instead consider the $C_4$-extension $K = \QQ(\alpha)$, with $\alpha^4 - \alpha^3 + \alpha^2 - \alpha + 1 = 0$,
then $K/\QQ$ is tamely ramified.

\begin{minted}{jlcon}
julia> K, a = number_field(x^4 - x^3 + x^2 - x + 1, "a");

julia> is_tamely_ramified(K)
true
\end{minted}

We do the same computation as before, but in this case obtain a free generator of the $\ZZ[G]$-module.

\begin{minted}{jlcon}
julia> V, f = galois_module(K); OK = ring_of_integers(K); M = f(OK);

julia> fl, c = is_free_with_basis(M); # the elements of c are a basis

julia> b = preimage(f, c[1]) # the element might different per session
a
\end{minted}

We can verify that the preimage of the element $c_1$ under the map $f$ is a normal integral basis generator.
As for the normal basis generator, we check whether the matrix of coordinates with respect to an integral basis is invertible over $\ZZ$:

\begin{minted}{jlcon}
julia> A, mA = automorphism_group(K);

julia> X = matrix(ZZ, [coordinates(OK(mA(s)(b))) for s in A])
[0    1    0    0]
[0    0    0    1]
[1   -1    1   -1]
[0    0   -1    0]

julia> det(X)
1
\end{minted}

One question we can ask is whether not being locally free is the only obstruction for a number field
to have a normal integral basis.
For abelian extensions of $\QQ$, the affirmative answer to this question is the content of the theorem of Hilbert--Speiser, which asserts
that an abelian extension $K/\QQ$ has a normal integral basis if and only if $K/\QQ$ is tamely ramified, or equivalently, there exists $n \in \NN$ odd and square-free such that $K \subseteq \QQ(\zeta_n)$.

For arbitrary extensions, the answer is in general a (complicated) no.
We next want to find a tamely ramified number field $K/\QQ$ without normal integral basis.
Since in this case $\mathcal{O}_K$ must be a locally free $\ZZ[G]$-module by Noether's criterion,
the obstruction to the existence of a normal integral basis must come from the existence
of locally free $\ZZ[G]$-modules which are not free.
An invariant capturing this obstruction (in part) is provided by the locally free class group $\Cl(\ZZ[G])$,
which is a finite abelian group and a generalization of the class group of Dedekind domains to the non commutative setting.
To have a chance of finding a locally free but non-free $\ZZ[G]$-module requires $\Cl(\ZZ[G]) \neq 1$.
Let us find the first non-abelian group where this is happening.

\begin{minted}{jlcon}
for i in 1:8, j in 1:number_of_small_groups(i)
  G = small_group(i, j)
  if is_abelian(G)
    continue
  end
  QG = QQ[G]
  ZG = integral_group_ring(QG)
  println(i, " ", j, " ", describe(G), ": ", locally_free_class_group(ZG))
end
\end{minted}

This yields:

\begin{minted}{jlcon}
6 1 S3: Z/1
8 3 D8: Z/1
8 4 Q8: Z/2
\end{minted}

We conclude that $G = \mathrm{Q}_8$ is the smallest group admitting locally free non-free $\ZZ[G]$-modules.
Next we construct a few candidate fields and hope that we find what we are looking for.
Since $\mathrm{Q}_8$ has the derived series $1 \subseteq \mathrm{C}_2 \subseteq \mathrm{Q}_8$ with quotient $\mathrm{Q}_8/\mathrm{C}_2 \cong \mathrm{C}_2 \times \mathrm{C}_2$, we can construct $\mathrm{Q}_8$-fields as quadratic extensions of $\mathrm{C}_2 \times \mathrm{C}_2$-fields.
As a starting field we take the $\mathrm{C}_2 \times \mathrm{C}_2$-extension $k/\QQ$ defined by $x^4 - 13x^2 + 16$.
On top of this field, using class field theory, we construct all tamely ramified $\mathrm{C}_2$-extensions $K/k$ which are normal over $\QQ$ and satisfy $\lvert \operatorname{disc}(K) \rvert \leq 10^{13}$ and $\operatorname{Gal}(K/\QQ) \cong \mathrm{Q}_8$.

\begin{minted}{jlcon}
k, = number_field(x^4 - 13*x^2 + 16)
candidates = []
for F in abelian_normal_extensions(k, [2], ZZ(10)^13)
  K, = absolute_simple_field(number_field(F))
  if !is_tamely_ramified(K)
    continue
  end
  if !is_isomorphic(galois_group(K)[1], quaternion_group(8))
    continue
  end
  push!(candidates, K)
end
\end{minted}

We can now investigate the second candidate:

\begin{minted}{jlcon}
julia> length(candidates)
2

julia> K = candidates[2]; V, f = galois_module(K); OK = ring_of_integers(K); M = f(OK);

julia> fl = is_free(M)
false

julia> defining_polynomial(K)
x^8 + 105*x^6 + 3465*x^4 + 44100*x^2 + 176400
\end{minted}

Thus we have found a tamely ramified $\mathrm{Q}_8$-extension $K/\QQ$ without normal integral basis.
In fact, with a little more effort we could prove that this is in fact the smallest $\mathrm{Q}_8$-extension with this property. For this we would have to consider all biquadratic fields $k$ up to the discriminant bound $10^{\frac{13}{2}}$.
For more information on the vast and intriguing subject of Galois module theory of rings of integers, we refer the reader to Fröhlich~\cite{MR717033}. Algorithmic aspects of this theory are treated for example in~\cite{MR2422318, MR2813368, MR4136552, MR4493243}.

\section*{Acknowledgements}
Funded by the Deutsche Forschungsgemeinschaft (DFG, German Research
Foundation) -- Project-ID 286237555 -- TRR 195. The first author also acknowledges funding through MaRDI (Mathematical Research
Data Initiative), funded by the Deutsche Forschungsgemeinschaft, project
number 460135501, NFDI 29/1 ``MaRDI -- Mathematische
Forschungsdateninitiative'' and the state of Rheinland-Pfalz via the
Forschungsinitiative and the ``SymbTools'' project.

\bibliographystyle{plain}
\bibliography{references}

 \end{document}